\documentclass[11pt]{article}
\usepackage{amsmath, amsfonts}

\numberwithin{equation}{section}

\def\B1{B_{1/2}}

\def\Bl1{{\lambda_1}}
\def\Box{\hfill\rule{2.5mm}{2.5mm}}

\def\H{{\cal H}}

\def\L{{\cal L}}

\def\N{{\mathbb {N}}}

\def\R{{\mathbb {R}}}
\def\RR{{\cal R}}
\def\SS{{\cal S}}

\def\am{{\alpha_-}}
\def\ap{{\alpha_1^i}}
\def\argth{\mathop {\rm argth}}

\def\bl1{{\bar \lambda_1}}

\def\build#1_#2^#3{\mathrel{
\mathop{\kern 0pt#1}\limits_{#2}^{#3}}}

\def\cd{{\bar d}}
\def\d{\displaystyle}

\def\E{{\cal E}}

\def\h1{\mathop{\rm H^1_{\rm loc,\rm u}}}

\def\iint{\displaystyle\int_{-1}^1}

\def\k{{m}}
\def\l2{\mathop{\rm L^2_{\rm loc,\rm u}}}

\def\lzero{{L_{\k+1}}}

\def\pr{{\partial_r}}
\def\ps{\partial_s}

\def\py{\partial_y}

\def\tmu{{\bar \nu}}

\def\ts{{s}}

\newcommand{\aref}[1]{(\ref{#1})}

\newcommand{\vc}[2]{
\left(
\begin{array}{l}
#1\\
#2
\end{array}
\right)
}

\title{\bf Blow-up behavior outside the origin for a semilinear wave equation in the radial case
}
\author{Frank Merle\\
{\it \small Universit\'e de Cergy Pontoise and IHES}\\
Hatem Zaag\\
{\it \small CNRS UMR 7539 LAGA Universit\'e Paris 13}}
\date{February 7, 2011}

\newtheorem{cor}{Corollary}[section]

\newtheorem{lem}[cor]{Lemma}
\newtheorem{prop}[cor]{Proposition}

\newtheorem{propo}{Proposition}
\newtheorem{theor}[propo]{Theorem}
\newtheorem{coro}[propo]{Corollary}

\begin{document}

\maketitle

{\small {\bf Abstract}: We consider the semilinear wave equation in the radial case with conformal subcritical power nonlinearity. If we consider a blow-up point different from the origin, then we exhibit a new Lyapunov functional which is a perturbation of the one dimensional case and extend all our previous results known in the one-dimensional case. In particular, we show that the blow-up set near non-zero non-characteristic points is of class $C^1$, and that the set of characteristic points is made of concentric spheres in finite number in $\{\frac 1R \le |x|\le R\}$ for any $R>1$.}

\medskip

{\bf MSC 2010 Classification}:  
35L05, 35L71, 35L67,
35B44, 35B40


\medskip

{\bf Keywords}: Wave equation, radial case, characteristic point, blow-up set.

\section{Introduction}
We consider radial solutions of the following semilinear wave equation:
\begin{equation}\label{equrn}
\left\{
\begin{array}{l}
\partial_t^2 U =\Delta U+|U|^{p-1}U,\\
U(0)=U_0\mbox{ and }U_t(0)=U_1,
\end{array}
\right.
\end{equation}
where $U(t):x\in\R^N \rightarrow U(x,t)\in\R$, $U_0\in \rm H^1_{\rm loc,u}$
and $U_1\in \rm L^2_{\rm loc,u}$.\\
 The space $\rm L^2_{\rm loc,u}$ is the set of all $v$ in 
$\rm L^2_{\rm loc}$ such that 
\[
\|v\|_{\rm L^2_{\rm loc,u}}\equiv\d\sup_{a\in \R^N}\left(\int_{|x-a|<1}|v(x)|^2dx\right)^{1/2}<+\infty,
\]
 and the space ${\rm H}^1_{\rm loc,u}= \{ v\;|\;v, \nabla v \in {\rm L}^2_{\rm loc,u}\}$.\\
 We assume in addition that 
\begin{equation}\label{condp}
1<p\mbox{ and }p\le 1+\frac 4{N-1}\mbox{ if } N\ge 2.
\end{equation}
Since $U$ is radial, we introduce
\begin{equation}\label{defu}
u(r,t) = U(x,t)\mbox { if }r=|x|
\end{equation}
and rewrite \eqref{equrn} as 
\begin{equation}\label{equ}
\left\{
\begin{array}{l}
\partial^2_{t} u =\partial^2_r u+\frac{(N-1)}r \pr u+|u|^{p-1}u,\\
\pr u(0,t)=0,\\
u(r,0)=u_0(r)\mbox{ and }u_t(r,0)=u_1(r),
\end{array}
\right.
\end{equation}
where $u(t):r\in\R^+ \rightarrow u(r,t)\in\R$.
We solve equation \eqref{equrn} locally in time in the space ${\rm H}^1_{\rm loc}\times {\rm L}^2_{\rm loc}(\R^N)$ (see Ginibre, Soffer and Velo \cite{GSVjfa92}, Lindblad and Sogge \cite{LSjfa95}).
Existence of blow-up solutions follows from ODE techniques or the energy-based blow-up criterion of Levine \cite{Ltams74}.
More blow-up results can be found in
Caffarelli and Friedman \cite{CFtams86}, \cite{CFarma85},
Alinhac \cite{Apndeta95}, \cite{Afle02} and Kichenassamy and Littman \cite{KL1cpde93}, \cite{KL2cpde93}.

\bigskip

If $u$ is a blow-up solution of \eqref{equ}, we define (see for example Alinhac \cite{Apndeta95}) a 1-Lipschitz curve $\Gamma=\{(r,T(r))\}$ where $r\ge 0$ 
such that the maximal influence domain $D$ of $u$ (or the domain of definition of $u$) is written as 
\begin{equation}\label{defdu}
D=\{(r,t)\;|\; t< T(r)\}.
\end{equation}
$\Gamma$ is called the blow-up graph of $u$. 
A point $r_0\ge 0$ is a non-characteristic point if there are 
\begin{equation}\label{nonchar}
\delta_0\in(0,1)\mbox{ and }t_0<T(r_0)\mbox{ such that }
u\;\;\mbox{is defined on }{\cal C}_{r_0, T(r_0), \delta_0}\cap \{t\ge t_0\}\cap\{r\ge 0\}
\end{equation}
where ${\cal C}_{\bar r, \bar t, \bar \delta}=\{(r,t)\;|\; t< \bar t-\bar \delta|r-\bar r|\}$. We denote by $\RR\subset \R^+$ (resp. $\SS\subset \R^+$) the set of non-characteristic (resp. characteristic) points.

\bigskip

In a series of papers \cite{MZjfa07}, \cite{MZcmp08}, \cite{MZajm10} and \cite{MZisol10} (see also the note \cite{MZxedp10}), we gave a full picture of the blow-up for solutions of \eqref{equrn} in one space dimension. In this paper, we aim at extending all those results to higher dimensions in the radial case, outside the origin. 
%
%
%
%
%
%
%

\bigskip

Throughout this paper, we consider $U(x,t)$ a radial blow-up solution of equation \eqref{equrn}, and use the notation $u(r,t)$ introduced in \eqref{defu}. We proceed in 3 sections:\\
- in Section \ref{seclyap}, we give a new Lyapunov functional for equation \eqref{equ} and bound the solution in the energy space.\\
- in Section \ref{secnonchar}, we study $\RR$, in particular the blow-up behavior of the solution and the regularity of the blow-up set there.\\
- in Section \ref{secchar}, we focus on $\SS$, both from the point of view of the blow-up behavior and the regularity of the blow-up set.

\section{A new Lyapunov functional and a new blow-up criterion in the radial case}\label{seclyap}
The extension of these results to higher dimensions in the radial setting involves the very beginning of our work in one dimension, namely the existence of a Lyapunov functional and the boundedness of the solution in similarity variables, performed in Antonini and Merle \cite{AMimrn01} and Merle and Zaag \cite{MZajm03}. 

\medskip

In \cite{HZlyap10} and \cite{HZlyapc10}, Hamza and Zaag considered the following class of perturbed wave equations for 
\begin{equation}\label{1} 
 u_{tt}= \Delta u+ |u|^{p-1} u+f(u)+g(u_t),\qquad  (x,t)\in \R^N\times \R_+^* 
\end{equation}
where $p$ satisfies \eqref{condp}
and for some $q\in[0,p)$,
\[
|{f(x)}|\le C(1+|x|^q) \mbox{ and }|{g(x)}|\le C(1+|x|).
\]
They showed in \cite{HZlyap10} and \cite{HZlyapc10} that the argument of Antonini, Merle and Zaag in \cite{AMimrn01}, \cite{MZajm03}, \cite{MZma05} and \cite{MZimrn05} extends through a perturbation method to equation \eqref 1. The key idea is to modify the Lyapunov functional of \cite{AMimrn01} with exponentially small terms and define a new functional which is in the same time decreasing in time and gives a blow-up criterion.

\medskip

It happens that the perturbation argument developed for equation \eqref 1  in \cite{HZlyap10} and \cite{HZlyapc10}  works for equation \eqref{equ} with no further modification, as far as blow-up points different from the origin are considered. For the reader's convenience, we briefly recall the argument of Hamza and Zaag in the following.

\bigskip

Given $r_0>0$, we recall the following similarity variables' transformation
\begin{equation}\label{defw}
w_{r_0}(y,s)=(T(r_0)-t)^{\frac 2{p-1}}u(r,t),\;\;y=\frac{r-r_0}{T(r_0)-t},\;\;
s=-\log(T(r_0)-t).
\end{equation}
The function $w=w_{r_0}$ satisfies the following equation for all $y\in (-1,1)$ and $s\ge \max\left(-\log T(r_0), -\log r_0\right)$:
\begin{equation}\label{eqw}
\partial^2_{s}w= \L w-\frac{2(p+1)}{(p-1)^2}w+|w|^{p-1}w
-\frac{p+3}{p-1}\partial_sw-2y\partial^2_{y,s} w+e^{-s}\frac{(N-1)}{r_0+ye^{-s}} \py w, 
\end{equation} 
\begin{equation}\label{defro}
\mbox{where }\L w = \frac 1\rho \py \left(\rho(1-y^2) \py w\right)\mbox{ and }
\rho(y)=(1-y^2)^{\frac 2{p-1}}.
\end{equation}
Let us recall the Lyapunov functional in one space dimension
\begin{equation}\label{defE}
E(w)= \iint \left(\frac 12 (\ps w)^2 + \frac 12  \left(\partial_y w\right)^2 (1-y^2)+\frac{(p+1)}{(p-1)^2}w^2 - \frac 1{p+1} |w|^{p+1}\right)\rho dy,
\end{equation}
which is defined in the space
\begin{equation}\label{defnh0}
\H = \left\{q \in {\rm H^1_{loc}}\times {\rm L^2_{loc}(-1,1)}
\;\;|\;\;\|q\|_{\H}^2\equiv \int_{-1}^1 \left(q_1^2+\left(q_1'\right)^2  (1-y^2)+q_2^2\right)\rho dy<+\infty\right\}.
\end{equation}
Introducing
\begin{equation}\label{defF}
F(w,s) = E(w) - e^{-s}\iint w \ps w  \rho dy,
\end{equation}
we claim the following:
\begin{prop}[A new functional for equation \eqref{eqw}]\label{prophamza}
$ $\\
(i) There exists $\gamma(p)>0$ and $S_0(N,p)\in\R$ such that for all $r_0>0$ and for all $s\ge \max\left(-\log T(r_0), S_0, S_0-2\log r_0\right)$,
\begin{equation}\label{edoF}
\frac d{ds}F(w_{r_0}(s),s)
\le \gamma(p)e^{-s}F(w_{r_0}(s),s)
-\frac 2{p-1}\iint (\ps w_{r_0}(s))^2 \frac{\rho}{1-y^2} dy.
\end{equation}
(ii) {\bf (A blow-up criterion)} There exists $S_1(p)\in \R$ such that if $W$ is a solution of equation (\ref{eqw}) with $\|W(s)\|_{L^{p+1}(B)}$ locally bounded, and $F(W(s_0),s_0)<0$ for some $s_0\ge S_1(p)$, then $W$ cannot be defined on the whole interval $[s_0, \infty)$. 
\end{prop}
{\bf Remark}: From (i), we see that the Lyapunov functional for equation \eqref{eqw} is in fact $H(w_{r_0}(s),s)$ where
\begin{equation}\label{defH}
H(w,s)=F(w,s)e^{-\gamma(p)e^{-s}},
\end{equation}
 not $F(w_{r_0}(s),s)$ nor $E(w_{r_0}(s))$.\\
{\bf Remark}: 
We already know from \cite{MZajm03} and \cite{MZma05} that even in the non-radial setting, equation \eqref{equ} has a Lyapunov functional given by a natural extension to higher dimensions of $E(w_{r_0}(s))$ \eqref{defE}. Unfortunately, due to the lack of information on stationary solutions in similarity variables in dimensions $N\ge 2$, we could not go further in our analysis, and we had to stop at the step of bounding the solution in similarity variables. On the contrary, when $N=1$, we could obtain a very precise characterization of blow-up.\\ 
 Here, restricting ourselves to radial solutions, we find a different Lyapunov functional in higher dimensions (which exists even for supercritical $p$). Considering blow-up points different from the origin, the characterization of stationary solutions in one space dimension is enough, and we are able to go in our analysis as far as in the one-dimensional case.

\bigskip

Following our analysis in \cite{MZajm03} and \cite{MZjfa07}, we derive with no difficulty the following:
\begin{prop}\label{boundedness}{\bf (Boundedness of the solutions of equation \eqref{eqw} in the energy space)} For all $r_0>0$, there is a $C_2(r_0)>0$ and $S_2(r_0)\in\R$ such that for all $r\in[\frac{r_0}2, \frac{3r_0}2]$ and $s\ge S_2(r_0)$,
\[
\iint\left((\py w_r(s))^2 (1-y^2) +  (w_r(s))^2
 + (\partial_s w_r(s))^2+|w_r(s)|^{p+1}\right) \rho dy\le C_2(r_0).
\]
 \end{prop}
{\it Proof}: The adaptation is straightforward from our analysis in \cite{MZajm03} and Proposition 3.5 page 66 in \cite{MZjfa07}. The only difference is in the justification of the limit at infinity of $E(w_{r_0}(s))$, which follows from the limit of $H(w_{r_0}(s),s)$ defined in \eqref{defH}. In fact, we know from Proposition \ref{prophamza} that $H(w_{r_0}(s),s)$ is decreasing and bounded from below, and such an information is unavailable for $E(w_{r_0}(s))$. 

\bigskip

{\it Proof of Proposition \ref{prophamza}}:\\
(i) Consider $r_0>0$, $s\ge \max(-\log T(r_0), -\log \frac{r_0}2)$ and write $w=w_{r_0}$ for simplicity. From the similarity variables' transformation \eqref{defw}, we see that 
\begin{equation}\label{range}
r=r_0+ye^{-s}\in\left[\frac{r_0}2, \frac{3r_0}2\right].
\end{equation}
Multiplying equation \eqref{eqw} by $\ps w\rho$ and integrating for $y\in(-1,1)$, we see by definition \eqref{defE} of $E(w)$ that
\begin{equation}\label{Eprime}
\frac d{ds}E(w(s)) =-\frac 4{p-1} \iint (\ps w(s))^2 \frac{\rho}{1-y^2} dy 
+(N-1) e^{-s} \iint \ps w(s)\py w(s) \frac{\rho}r dy
\end{equation}
where $r$ is defined in \eqref{range}.
Using \eqref{range}, we write
\begin{eqnarray}
&&\left|(N-1)  \iint \ps w(s)\py w(s) \frac{\rho}r dy\right|\label{CS}\\
&\le& \left(\frac{p-1}4\right)\iint (\py w(s))^2\rho(1-y^2) dy +  
\left(\frac 4{p-1}\right)\left(\frac{N-1}{r_0}\right)^2\iint (\ps w(s))^2 \frac{\rho}{1-y^2} dy.\nonumber
\end{eqnarray}
Recalling the following Hardy-Sobolev estimate (see Appendix B page 1163 in \cite{MZajm03} for the proof):
\begin{equation}\label{hs}
\iint h^2 \frac \rho{1-y^2} dy \le C\iint h^2 \rho dy + C \iint (h'(y))^2 \rho(1-y^2)dy,
\end{equation}
we use the notation $I(s) = e^{-s}\iint \ps w w \rho dy$ and write from equation \eqref{eqw} for any $\epsilon>0$ and $s\ge -\log\left(\frac{\epsilon r_0}{2(N-1)}\right)$,
\begin{eqnarray}
e^s I'(s) &\ge& \iint |w(s)|^{p+1} \rho dy -(1+\epsilon)\iint |\py w(s)|^2\rho(1-y^2) dy\label{Iprime}\\
&& -\frac {C(p)}\epsilon \iint |\ps w(s)|^2 \frac{\rho}{1-y^2} dy-(\epsilon +\frac{2(p+1)}{(p-1)^2})\iint (w(s))^2 \rho dy.\nonumber
\end{eqnarray}
Using \eqref{Eprime}, \eqref{CS} and \eqref{Iprime}, we get that \eqref{edoF} follows by definition \eqref{defF} of $F(w,s)$, 
provided that we fix $\epsilon=\epsilon(p)>0$ small enough and take $s$ large enough. This yields (i) of Proposition \ref{prophamza}.

\medskip

(ii) If $W$ is a solution of equation \eqref{eqw}, then we write by definition \eqref{defH} of $H$:
\begin{eqnarray*}
&&H(W(s),s)\ge -\frac{e^{-\gamma(p)e^{-s}}}{p+1}\iint |W(s)|^{p+1} \rho dy\\
&& +e^{-\gamma(p)e^{-s}}\left(\left(\frac 12 - \frac{e^{-s}}2\right)\iint (\ps W(s))^2 \rho dy +\left(\frac{p+1}{(p-1)^2} - \frac{e^{-s}}2\right)\iint (W(s))^2 \rho dy\right)\\
&\ge & -\frac 1{p+1}\iint |W(s)|^{p+1} \rho dy
\end{eqnarray*} 
if $s\ge S_1(p)$ for some $S_1(p)\in\R$ large enough. Using this inequality together with the fact that $H(W(s),s)$ is decreasing by the remark following Proposition \ref{prophamza}, we see that the argument used by Antonini and Merle in Theorem 2 page 1147 in \cite{AMimrn01} for the equation \eqref{equrn} works here and we get the blow-up criterion. This concludes the proof of Proposition \ref{prophamza}.\Box
\section{Blow-up results related to non-characteristic points}\label{secnonchar}
Let us first introduce for all $|d|<1$ the following solitons defined by 
\begin{equation}\label{defkd}
\kappa(d,y)=\kappa_0 \frac{(1-d^2)^{\frac 1{p-1}}}{(1+dy)^{\frac 2{p-1}}}\mbox{ where }\kappa_0 = \left(\frac{2(p+1)}{(p-1)^2}\right)^{\frac 1{p-1}} \mbox{ and }|y|<1.
\end{equation}
Note that $\kappa(d)$ is a stationary solution of \eqref{eqw} in one space dimension. 

\medskip

\noindent Adapting the analysis of \cite{MZjfa07} and \cite{MZcmp08}, we claim the following:
\begin{theor}[Blow-up behavior and regularity of the blow-up set on $\RR$]\label{thbb}$ $\\
(i) {\bf (Regularity related to $\RR$)} $\RR\neq \emptyset$, $\RR\cap \R^*_+$ is an open set, and $x\mapsto T(x)$ is of class $C^1$ on $\RR\cap \R^*_+$.\\
(ii) {\bf (Blow-up behavior in similarity variables)} There exist $\mu_0>0$ and $C_0>0$ such that for all $r_0\in\RR\cap \R^*_+$, there exist
$\theta(r_0)=\pm 1$ and $s_0(r_0)\ge - \log T(r_0)$ such that for all $s\ge s_0$:
\begin{equation}\label{profile}
\left\|\vc{w_{r_0}(s)}{\partial_s w_{r_0}(s)}-\theta(r_0)\vc{\kappa(T'(r_0))}{0}\right\|_{\H}\le C_0 e^{-\mu_0(s-s_0)}.
\end{equation}
Moreover, $E(w_{r_0}(s)) \to E(\kappa_0)$ as $s\to \infty$.
\end{theor}
{\bf Remark}: If $0\in\RR$, the asymptotic behavior of $w_0$ remains open.\\
{\it Proof}: The proof is the same as for the one-dimensional case treated in \cite{MZjfa07} and \cite{MZcmp08}, with some minor adaptations. For that reason, we don't give the full proof here. We will instead ask the reader to follow our proof in those papers, and provide him only with the delicate points he may face in adapting the proof to the higher-dimensional radial case.\\
As in the one-dimensional case, we proceed in two steps: we first prove (ii) with a parameter $d_0(r_0)$ instead of $T'(r_0)$, then we prove (i) and the fact that in (ii) $d_0(r_0)=T'(r_0)$.

\bigskip

{\bf Step 1: Proof of (ii) with $d_0(r_0)$ instead of $T'(r_0)$ in \eqref{profile}}

This statement is the twin of Corollary 4 page 49 given in \cite{MZjfa07} in the one-dimensional case.\\
Consider $r_0\in\RR \cap \R^*_+$. 
As in that case, the proof has two major steps:

\medskip

- {\it Approaching the set of stationary solutions}, in the sense that for some $\theta(r_0)=\pm 1$ and $d_0(r_0)\in(0,1)$, we have
\begin{equation}\label{appr}
\inf_{|d|\le d_0(r_0)}\|w_{r_0}(\cdot,s)- \theta(r_0) \kappa(d,\cdot)\|_{H^1(-1,1)}+\|\partial_s w_{x_0}\|_{L^2(-1,1)}\to 0
\end{equation}
 as $s\to \infty$ (see Theorem 2 page 47 of \cite{MZjfa07} for the statement in 
 one space dimension). Note that such a statement is still unavailable for non radial solutions of equation \eqref{equrn}, since we have no classification for stationary solutions of the multi-dimensional version of equation \eqref{eqw}. Fortunately, in the radial case, we naturally see from equation \eqref{eqw} and the estimate \eqref{range} that we only need the classification in one space dimension, given in Proposition 1 page 46 in \cite{MZjfa07}. This is the reason why we restrict ourselves to the radial case in this paper. As for the proof of \eqref{appr}, the reader has to see Section 3 page 60 in \cite{MZjfa07}. The only delicate point is the adaptation of the proofs of Lemma 3.3 page 62 and Lemma 3.4 page 64. Indeed:\\
- In the proof of Lemma 3.3, we need a Duhamel formulation for the one-dimensional semilinear wave equation (see estimate (63) page 63 in \cite{MZjfa07}). This formulation has to be replaced by the radial version, which may be derived from Shatah and Struwe \cite{SSnyu98}.\\
%
%
- As for the proof of Lemma 3.4, due to the fact that equation \eqref{equ} is no longer invariant under the scaling
\[
\lambda \mapsto u_\lambda(\xi, \tau) = \lambda^{\frac 2{p-1}}u(\lambda \xi, \lambda \tau),
\]
we need to understand the continuous dependence of the solutions of the following family of equations 
\begin{equation}\label{eqxl}
\partial^2_{t} u =\partial^2_r u+\frac{(N-1)\lambda}{x+\lambda r} \pr u+|u|^{p-1}u,
\end{equation}
with respect to initial data and the parameters $x\ge 0$ and $\lambda>0$ (including the limit as $\lambda \to 0$), and this is a classical estimate. Apart from this point, there is no other problem in the adaptation. 

\medskip 

- {\it Proof of the convergence to a $\theta(r_0) \kappa(d_0(r_0))$ for some $d_0(r_0)\in(-1,1)$}. From \eqref{appr} and the monotonicity of the Lyapunov function $H(w_{r_0}(s),s)$ \eqref{defH}, we see that (ii) (with a parameter $d_0(r_0)$ instead of $T'(r_0)$ in \eqref{profile}) follows directly from the following trapping result. Note that the statement of this result is different from the analogous result in one space dimension given in Theorem 3 page 48 in \cite{MZjfa07}. This difference is due to the fact that the equation \eqref{eqw} in similarity variables depends on a parameter $r_0>0$ and contains a new term of order $e^{-s}$ (it is no longer autonomous). 
For a further purpose, we give in the following the radial case's version of the trapping result, valid uniformly for all $r_0\ge \rho_0$ for some $\rho_0>0$:
\begin{theor}{\bf (Trapping near the set of non zero stationary solutions of \eqref{eqw})}\label{thtrap}
\\
For all $\rho_0>0$, there exist positive $\epsilon_0$, $\mu_0$ and $C_0$ such that for all $\epsilon^*\le \epsilon_0$, there exists $s_0(\epsilon^*)$ such that if $r_0\ge \rho_0$, $s^*\ge s_0$ and $w\in C([s^*, \infty), \H)$ is a solution of equation \eqref{eqw} with
\begin{equation}\label{highenergy}
\forall s\ge s^*,\;\;E(w(s))\ge E(\kappa_0)-e^{-\frac s2},
\end{equation}
and 
\begin{equation*}
\left\|\vc{w(s^*)}{\partial_s w(s^*)}-\omega^*\vc{\kappa(d^*,\cdot)}{0}\right\|_{\H}\le \epsilon^*
\end{equation*}
for some $d^*\in(-1,1)$ and $\omega^*=\pm 1$,
then there exists $d_\infty\in (-1, 1)$ such that
\[
\left|\argth{d_\infty} - \argth{d^*}\right|\le C_0 \epsilon^*,
\]
and for all $s\ge s^*$,
\begin{equation*}
\left\|\vc{w(s)}{\partial_s w(s)}-\omega^*\vc{\kappa(d_\infty, \cdot)}{0}\right\|_{\H}\le C_0 \epsilon^* e^{-\mu_0(s-s^*)}.
\end{equation*}
\end{theor}
{\bf Remark}: The function $\argth$ is the inverse of the hyperbolic tangent function. It is given by $\argth d = \frac 12 \log\left(\frac{1+d}{1-d}\right)$. Theorem \ref{thtrap} holds under a weaker condition where we replace the $e^{-\frac s2}$ in \eqref{highenergy} by $C^* e^{-s}$ for some constant $C^*(\rho_0)>0$.\\
{\it Proof}: The proof can be adapted from the proof of Theorem 3 in \cite{MZjfa07} given in Section 5 page 99 in that paper. Up to replacing $u$ by $-u$, we may assume $\omega^*=1$. As in \cite{MZjfa07}, we linearize equation \eqref{eqw} by introducing
\begin{equation}\label{defq0}
q(y,s)=w(y,s) - \kappa(d(s),y)
\end{equation}
where the parameter $d(s) \in(-1,1)$ is chosen by modulation so that
\begin{equation}\label{kill}
\pi_0^{d(s)}(q(s))=0,
\end{equation}
and $\pi_0^d$ is the projector on the null mode of the linearized operator. Then, we decompose $q$ into two parts, according to the spectrum of the linearized operator as follows:\\
- its projection on the mode $\lambda=1$, whose norm is bounded by some $|\alpha_1(s)|$;\\
- its projection on the negative part of the spectrum, whose norm is bounded by some $\alpha_-(s)$.\\
Thanks to this decomposition, we write
\[
\frac 1{C_0}\left(|\alpha_1(s)|+\alpha_-(s)\right)\le \|q(s)\|_{\H}\le C_0\left(|\alpha_1(s)|+\alpha_-(s)\right)
\]
for some $C_0>0$.\\
Projecting the linearized equation according to this decomposition, we see that we have exponentially small perturbations coming from the new term in the similarity variables equation \eqref{eqw}, which is uniformly bounded since the parameter $r_0$ satisfies $r_0\ge \rho_0>0$, in the sense that
\begin{equation}\label{small}
\left|\frac{(N-1)e^{-s}}{r_0+ye^{-s}}\py w\right|\le \frac 2{\rho_0} (N-1) e^{-s}|\py w|\mbox{ as soon as } s\ge -\log \frac{\rho_0}2.
\end{equation}
More precisely, we have the following statement, which is the new version of Proposition 5.2 page 103 in \cite{MZjfa07}:
\begin{prop}\label{lemproj}
There exists $\epsilon_2>0$ such that if $w$ a solution to equation \eqref{eqw} satisfying \eqref{kill} and 
\[
\|q(s)\|_{\H} \le \epsilon
\]
at some time $s$ for some $\epsilon\le \epsilon_2$, where $q$ is defined in \eqref{defq0}, then:\\
(i) {\bf (Control of the modulation parameter)}
\begin{equation*}
|d'|\le C_0(1-d^2)(\ap^2+\am^2)+Ce^{-s}(1-d^2).
\end{equation*}
(ii) {\bf (Projection of the linearized equation on the different eigenspaces of $L_d$)} 
\begin{eqnarray*}
\left|\ap'- \ap\right|&\le& C_0\left(\ap^2+\am^2\right)+Ce^{-s},
\\
\left(R_-+\frac 12 \am^2\right)'&\le &- \frac 7{2(p-1)}\iint q_{-,2}^2 \frac \rho{1-y^2} dy
+ C_0\left(\ap^2+\am^2\right)^{3/2}+Ce^{-s}
\end{eqnarray*}
for some $R_-(s)$ satisfying
\begin{equation*}
|R_-(s)|\le C_0 (\ap^2+\am^2)^{\frac{1+\bar p}2}\mbox{ where } \bar p= \min(p,2) >1.
\end{equation*}
(iii) {\bf (Additional relation)}
\begin{equation*}
\frac d{ds}\iint q_1q_2 \rho \le 
-\frac 45 \am^2+C_0 \iint q_{-,2}^2 \frac \rho{1-y^2}+C_0\ap^2+Ce^{-s}.
\end{equation*}
(iv) {\bf (Energy barrier)} If moreover \eqref{highenergy} holds, then
\begin{equation*}
\ap(s)\le C_0 \am(s)+Ce^{-s/2}.
\end{equation*}
\end{prop}
{\it Proof}: In comparison with the one dimensional case, the linearized equation has an additional exponentially small term. Therefore, the adaptation consists in projecting that term on the different components, and there is no difficulty for this. The projection of the remaining terms is the same as the one-dimensional case. See the ``Proof of Proposition 5.2'' page 104 in \cite{MZjfa07} for the one-dimensional case. This concludes the instructions for the adaptation of Proposition \ref{lemproj}.\Box

\medskip

Now, with this proposition, there is no particular difficulty in deriving Theorem \ref{thtrap} as in the one-dimensional case. See Section 5.3 page 113 in \cite{MZjfa07}. This concludes the instructions for the adaptation of Theorem \ref{thtrap}. \Box 

\bigskip

{\bf Step 2: Proof of  (i) and the fact that (ii) holds with a parameter $d_0(r_0)=T'(r_0)$ in \eqref{profile}}

Note that the fact that $\RR\neq \emptyset$ follows by the same argument given in the remark following Theorem 1 page 58 in \cite{MZcmp08}.\\
 The corresponding statement to this step in the one-dimensional case is Theorem 1 page 58 in \cite{MZcmp08}. As in that case, we need a rigidity Theorem (or a Liouville Theorem) in similarity variables. Let us insist on the fact that we don't need any adaptation for the Liouville Theorem and that the version we need is indeed the one-dimensional version stated in Theorems 2 and 2' in pages 58 and 59 of \cite{MZcmp08}. 
While adapting the proofs to the radial case, the reader should pay attention to two facts:\\
- Again, we need to use the continuity of solutions of equation \eqref{eqxl} with respect to initial data and to the parameters $(x,\lambda)$.\\
- The functional $E(w)$ \eqref{defE} is no longer nondecreasing. Fortunately, we can replace it by the Lyapunov functional $H(w,s)$ \eqref{defH} which is decreasing. The error we make in this replacement is exponentially small as one sees from \eqref{defF}, \eqref{defH} and Proposition \ref{boundedness}.\\
- We need the trapping result here too. Of course, the new version given in Theorem \ref{thtrap} applies.

\medskip

 Apart from these three remarks, the adaptation is straightforward. One has just to follow the proof of Theorem 1 given in Section 2 page 60 in \cite{MZcmp08}.

\section{Blow-up results related to characteristic points}\label{secchar}
In Proposition 1 in \cite{MZajm10}, we showed the existence of a solution to equation \eqref{equrn} in one space dimension such that $\SS \neq \emptyset$. Artificially adding other coordinates, this one-dimensional solution can be considered as a multi-dimensional solution with $\SS \neq \emptyset$. Adapting the argument of \cite{MZajm10} to the radial case, we prove the following:
\begin{propo}\label{propexis}{\bf (Existence of radial solutions with a non zero characteristic point)} There exists $(u_0,u_1)$ such that the corresponding solution of equation \eqref{equ} has a non zero characteristic point.
\end{propo}
{\bf Remark}: In this case, the multi-dimensional version $U(x,t)=u(|x|,t)$ has a sphere of characteristic points.\\
In fact, Proposition \ref{propexis} is a consequence of the following:
\begin{theor}[Existence and generic stability of characteristic points]
\label{thexis}
$ $\\
(i) {\bf (Existence)} Let $0<a_1<a_2$ be two non characteristic points such that
\[
w_{a_i}(s) \to \theta(a_i)\kappa(d_{a_i},\cdot)\mbox{ as }s\to \infty\mbox{ with }\theta(a_1)\theta(a_2)=-1
\]
for some $d_{a_i}$ in $(-1,1)$, in the sense \eqref{profile}. Then, there exists a characteristic point $c\in (a_1,a_2)$.\\
(ii) {\bf (Stability)} There exists $\epsilon_0>0$ such that if $\|(\tilde U_0,\tilde U_1)- (U_0, U_1)\|_{\h1\times \l2(\R^N)}\le \epsilon_0$, then, $\tilde u(r,t)$ the solution of equation \eqref{equ} with initial data $(\tilde u_0,\tilde u_1)(r)=(\tilde U_0,\tilde U_1)(x)$ if $r=|x|$ blows up and has a characteristic point $\tilde c\in [a_1,a_2]$. 
\end{theor}
{\bf Remark}: This statement is different from the original one (Theorem 2 in \cite{MZajm10}) by two natural small facts:
we take positive points $a_1$ and $a_2$ in (i), and 
we use the multi-dimensional norm in (ii) (of course, from the finite speed of propagation, it is enough to take a localized norm instead).\\
Indeed, let us first derive Proposition \ref{propexis} from this theorem, then prove this latter.

\medskip

{\it Proof of Proposition \ref{propexis} assuming Theorem \ref{thexis}}: From the finite speed of propagation and the uniqueness of the solution to the Cauchy problem, one can take $(u_0,u_1)$ with large plateaus of opposite signs so that for some $0<a_1<a_2$ and large $T>0$, the solution remains space independent in the backward cones $\{r-a_i<T-t\}$, conserves its sign there, and blows up at $a_i$ at some time $T(a_i)<T$.
Since the solution is space independent around $a_i$, the blow-up time is locally constant and the point $a_i\in\RR$ with $T'(a_i)=0$. Using the description of the blow-up behavior in the non-characteristic case stated in (ii) of Theorem \ref{thbb}, we see that $w_{a_i}(y,s) \to\theta(a_i) \kappa$ as $s\to \infty$ with $\theta(a_1) \theta(a_2) = -1$, since the plateaus have opposite signs. Therefore, the hypotheses of (i) of Theorem \ref{thexis} are fulfilled and we have the desired conclusion of Proposition \ref{propexis}.\Box

\bigskip

Now, we give indications on the adaptation of the proof of Theorem \ref{thexis} from the one-dimensional case.

\medskip

{\it Proof of Theorem \ref{thexis}}: There is no difficulty in adapting to the present context the proof of Theorem 2 of \cite{MZajm10} given in Section 2 of that paper, except may be for the continuity of the blow-up time with respect to initial data, stated in Proposition 2.1 of \cite{MZajm10}, where some natural extensions to the radial case are needed.\Box 

\bigskip

We also have the following result which relates the existence of characteristic points to the sign-change of the solution:
\begin{theor}
{\bf (Non-existence of characteristic points if the sign is constant)} Consider $u(r,t)$ a blow-up solution of \eqref{equ} such that $u(r,t)\ge 0$ for all $r\in (a_0,b_0)$ and $t_0\le t< T(r)$ for some real $0\le a_0<b_0$ and $t_0\ge 0$. Then, $(a_0, b_0)\subset \RR$.
\end{theor}
{\bf Remark}: This statement is exactly the same as the original (Theorem 4 in \cite{MZajm10}). In particular, it is valid with $a_0=0$.\\
{\it Proof}: This result follows from Theorem \ref{new} below exactly as in one space dimension. See the proof of Theorem 4 given in Section 4.1 in \cite{MZajm10}.\Box

\bigskip

Now, given $r_0\in\SS \cap \R^*_+$, we have the same description for the asymptotics of $w_{r_0}$ as in the one-dimensional case. More precisely, the following holds in the radial case, outside the origin (for the statement in one space dimension, see Theorem 6 in \cite{MZajm10}): 
\begin{theor}\label{new}
{\bf (Description of the behavior of $w_{r_0}$ where $r_0$ is characteristic)} Consider $u(r,t)$ a blow-up solution of \eqref{equ} and $r_0\in \SS\cap \R^*_+$. Then,  it holds that
\begin{equation}\label{cprofile00}
\left\|\vc{w_{r_0}(s)}{\ps w_{r_0}(s)} - \vc{\d\sum_{i=1}^{k(r_0)} e^*_i\kappa(d_i(s),\cdot)}0\right\|_{\H} \to 0\mbox{ and }E(w_{r_0}(s))\to k(r_0)E(\kappa_0)
\end{equation}
as $s\to \infty$, for some 
\begin{eqnarray}
k(r_0)&\ge& 2,\label{pb}\\
e^*_i&=&e^*_1(-1)^{i+1}\label{ei}
\end{eqnarray}
 and continuous $d_i(s)=-\tanh \zeta_i(s)\in (-1,1)$ for $i=1,...,k(r_0)$. Moreover, for some $C_0>0$, for all $i=1,...,k(r_0)$ and $s$ large enough, we have
\begin{equation}\label{equid}
\left|\zeta_i(s)-\left(i-\frac{(k(r_0)+1)}2\right)\frac{(p-1)}2\log s\right|\le C_0. 
\end{equation}
\end{theor}
{\it Proof}: As in the one-dimensional case, the proof of the asymptotic behavior and the geometric results on $\SS$ (see Theorem \ref{thgeo} below) go side by side. For that reason, we leave the proof after the statement of Theorem \ref{thgeo}.\Box



\bigskip

Extending the definition of $k(r_0)$ defined on $\SS$ in Theorem \ref{new} by setting
\[
\forall r_0\in\RR,\;\;k(r_0)=1,
\]
we get the following result on the energy behavior from the asymptotic behavior at a non-characteristic point (see (ii) of Theorem \ref{thbb}) and at a characteristic point (see Theorem \ref{new}):
\begin{coro}[A criterion for non characteristic points]\label{corcriterion}$ $\\
 For all $r_0>0$, there exist $C_3(r_0)>0$ and $S_3(r_0)\in\R$ such that:\\
(i) For all $r\in[\frac{r_0}2, \frac{3r_0}2]$ and $s\ge S_3$, we have
\[
E(w_r(s))\ge k(r)E(\kappa_0)-C_3(r_0)e^{-s}.
\]
(ii) If for some $r\in[\frac{r_0}2, \frac{3r_0}2]$ and $s\ge S_3$, we have 
\[
E(w_r(s))<2 E(\kappa_0)-C_3(r_0)e^{-s},
\]
then $r\in \RR$.
\end{coro}
{\bf Remark}: With respect to the one-space dimension statement (Corollary 7 in \cite{MZajm10}), this statement has additional exponentially small terms. This comes from the fact that the functional $E(w)$ is no longer decreasing, and that one has to work instead with the functional $H(w,s)$ \eqref{defH} which is decreasing, and differs from $E(w)$ by exponentially small terms, uniformly controlled for $r\in[\frac{r_0}2, \frac{3r_0}2]$ thanks to the uniform estimates of Proposition \ref{boundedness}.\\
{\it Proof}: If one replaces $E(w)$ by $H(w,s)$, then the proof is straightforward from Theorems \ref{thbb} and \ref{new} together with the monotonicity of $H(w,s)$ (see \eqref{defH} and \eqref{prophamza}). Since the difference between the two functionals is exponentially small, uniformly for $r\in[\frac{r_0}2, \frac{3r_0}2]$ (see \eqref{defH}, \eqref{edoF} and Proposition \ref{boundedness}), we get the conclusion of Corollary \ref{corcriterion}.\Box 

\bigskip

Finally, we give in the following some geometric information related to characteristic points (for the statement in one space dimension, see Theorem 1, Theorem 2 and the following remark in \cite{MZisol10}):
\begin{theor}\label{thgeo}{\bf (Geometric considerations on $\SS$)} Consider $u(r,t)$ a blow-up solution of equation \eqref{equ}.\\
(i) {\bf (Isolatedness of characteristic points)} All characteristic points different from the origin are isolated.\\
(ii) {\bf (Corner shape of the blow-up curve at characteristic points)} If $r_0\in \SS\cap \R^*_+$ and $0<|r-r_0|\le \delta_0$, then 
\begin{equation*}
\d\frac{1}{C_0|\log(r-r_0)|^{\frac{(k(r_0)-1)(p-1)}2}}\le T'(r)+\frac{r-r_0}{|r-r_0|} \le \frac{C_0}{|\log(r-r_0)|^{\frac{(k(r_0)-1)(p-1)}2}}
\end{equation*}
for some $\delta_0>0$ and $C_0>0$, where $k(r_0)\ge 2$ is the integer defined in \eqref{pb}.
\end{theor}
{\it Proof}: See below.\\
{\bf Remark}: Integrating the estimate in (ii) of this theorem, we see that
\begin{equation}\label{corner}
\d\frac{|x-x_0|}{C_0|\log(x-x_0)|^{\frac{(k(x_0)-1)(p-1)}2}}\le T(x)- T(x_0)+|x-x_0| \le \frac{C_0|x-x_0|}{|\log(x-x_0)|^{\frac{(k(x_0)-1)(p-1)}2}}.
\end{equation}
{\bf Remark}: Note from (i) that the multi-dimensional version $U(x,t)=u(|x|,t)$ has a finite number of concentric spheres of characteristic points in the set $\{\frac 1R < |x|<R\}$ for every $R>1$. This is consistent with our conjecture in \cite{MZisol10} where we guessed that in dimension $N\ge 2$, the $(N-1)$-dimensional Hausdorff measure of $\SS$ is bounded in compact sets of $\R^N$. Note that this conjecture is related to the result of Vel\'azquez who proved in \cite{Viumj93} that the $(N-1)$-dimensional Hausdorff measure of the blow-up set for the semilinear heat equation with subcritical power nonlinearity is bounded in compact sets of $\R^N$.

\bigskip

As a consequence of our analysis, particularly the lower bound on $T(r)$ in \aref{corner}, we have the following estimate on the blow-up speed in the backward light cone with vertex $(r_0, T(r_0))$ where $r_0> 0$ (for the statement in one space dimension, see Corollary 3 in \cite{MZisol10}): 
\begin{coro}\label{corspeed}{\bf (Blow-up speed in the backward light cone)} For all $r_0>0$, there exists $C_4(r_0)>0$ such that for all $t\in[0, T(r_0))$, we have
\[
\frac{|\log(T(r_0)-t)|^{\frac{k(r_0)-1}2}}{C_4(r_0)(T(r_0)-t)^{\frac 2{p-1}}}\le \sup_{|x-r_0|<T(r_0)-t}|u(x,t)|\le \frac{C_4(r_0) |\log(T(r_0)-t)|^{\frac{k(r_0)-1}2}}{(T(r_0)-t)^{\frac 2{p-1}}}.
\]
\end{coro}
{\bf Remark}: Note that when $r_0\in\RR\cap \R^*_+$, the blow-up rate of $u$ in the backward light cone with vertex $(r_0, T(r_0))$ is given by the solution of the associated ODE $u"=u^p$. When $r_0\in\SS\cap \R^*_+$, the blow-up rate is higher and quantified, according to $k(r_0)$, the number of solitons appearing in the decomposition \eqref{cprofile00}.\\
{\it Proof}: When $r_0\in\RR$, the result follows from the fact that the convergence in \aref{profile} is true also in $L^\infty\times L^2$ from \aref{appr} and the Sobolev embedding in one dimension. When $r_0\in\SS$, see the proof of Corollary 3 of \cite{MZisol10} given in Section 3.3 of that paper.\Box

\bigskip

{\it Proof of Theorems \ref{new} and \ref{thgeo}}: The proof follows the pattern of the original proof, given in \cite{MZjfa07}, \cite{MZajm10} and \cite{MZisol10}. In the following, we recall its different parts.

\bigskip

 {\bf Part 1: Proof of \eqref{cprofile00} without \eqref{pb} nor \eqref{ei} and with the estimate 
\begin{equation}\label{decuple}
\zeta_{i+1}(s)-\zeta_i(s) \to \infty\mbox{ as }s\to\infty
\end{equation}
 instead of \eqref{equid}} (note that both \eqref{equid} and \eqref{decuple} are meaningful only when $k(r_0)\ge 2$).

 The original statement of this part is given in Theorem 2 (B) page 47 in \cite{MZjfa07} and the proof in section 3.2 page 66 in that paper. Note that this part doesn't exclude the possibility of having $k(r_0)=0$ or $k(r_0)=1$. The adaptation is straightforward. As in the non-characteristic case above, one has to use the Duhamel formulation in the radial which may be derived from \cite{SSnyu98}.


\bigskip

{\bf Part 2: Proof of  \eqref{ei}, \eqref{equid} and the upper bound in \aref{corner}, assuming that \eqref{pb} is true}.

The original statement is given in Propositions 3.1 and 3.13 in \cite{MZajm10}. The reader has to read Section 3 and Appendices B and C in that paper.
The adaptation is straightforward, except for the effect of the new term in equation \eqref{eqw}, which produces exponentially small terms in many parts of the proof. In particular, Lemma 3.11 of \cite{MZajm10} has to be changed by adding $Ce^{-s}$ to the right of all the differential inequalities.

\bigskip

 {\bf Part 3: Proof of \eqref{pb} and the fact that the interior of $\SS$ is empty}. 

The original statement is given in Proposition 4.1 of \cite{MZajm10}. Here, the adaptation is not only delicate, but we need a new argument to rule out the occurrence of the case where, locally near the origin, the blow-up set of the multi-dimensional version $U(x,t)$ is a forward light cone with vertex $(0,T(0))$ (we actually prove a stronger result, see Lemma \ref{newlem} below). For this reason, we will give the good version of Proposition 4.1 in \cite{MZajm10} and outline its proof in the radial case. More precisely, we claim the following:
\begin{prop}\label{pbut}$ $\\
(i) The interior of $\SS$ is empty.\\
(ii) For all $r_0\in \SS\cap \R^*_+$, $k(r_0)\ge 2$.
\end{prop}
{\bf Remark}: Please note that in (i), the information is about $\SS$, whereas in (ii), we have to restrict to $\SS\cap \R^*_+$.

\medskip

As in \cite{MZajm10}, this proposition is a consequence of the following Lemmas, which we restate, since some statements surprisingly remain valid even at the origin, whereas others are valid only outside the origin:
\begin{lem}\label{cns}{\bf (Characterization of the interior of $\SS$)}
For any $0\le r_1<r_2$, the following statements are equivalent:

(a) $(r_1, r_2) \in \SS$.

(b) There exists $r^*\in [r_1, r_2]$ such that for all $r\in [r_1, r_2]$, $T(r) = T(r^*) - |r-r^*|$.
\end{lem}
{\it Remark and Proof}: Note that this lemma is valid also at the origin. The proof is the same as for Lemma 4.2 in \cite{MZajm10}.\Box
\begin{lem}\label{maxslope}
Consider $0\le r_1<r_2$ such that $e\equiv \frac{T(r_2)-T(r_1)}{r_2-r_1}=\pm 1$. Then,\\
(i) for all $r\in [r_1, r_2]$, $T(r)=T(r_1)+e(r-r_1)$,\\
(ii) $(r_1, r_2) \in \SS$.
\end{lem}
{\it Remark and Proof}: Note that this lemma is valid at the origin too. The proof can be adapted straightforwardly from the proof of Lemma 4.3 in \cite{MZajm10}.\Box 
\begin{lem}[Boundary properties of $\SS$]\label{lemfrank}$ $\\
(i) For all $r_0\in \partial \SS\cap \R^*_+$, $k(r_0) \neq 0$.\\
(ii) Consider $r_0\in \partial \SS\cap \R^*_+$ with $k(r_0)=1$.
If there exists a sequence $r_n\in \RR$ converging from the left (resp. the right) to $r_0$, then $r_0$ is left-non-characteristic (resp. right-non-characteristic).
\end{lem}
{\bf Remark}: Unlike the two previous lemmas, this lemma is valid outside the origin. This is due to the fact that we strongly need the structure in similarity variables, which is available only outside the origin. We mean by $r_0$ is left-non-characteristic (resp. right-non-characteristic) that it satisfies condition \eqref{nonchar} only for $r<r_0$ (resp. for $r>r_0$).\\
{\it Proof}: The only delicate point in the adaptation of the proof from the proof of Lemma 4.4 in \cite{MZajm10} is in Claim 4.5. Indeed, we need to choose there the time $\tilde t$ close enough to $T(x_0)$ so that we can apply the trapping result stated in Theorem \ref{thtrap}. Remember that this restriction in the trapping result comes from the fact that the equation \eqref{eqw} is no longer autonomous, and that the new term in \eqref{eqw} becomes small when $s$ is large (see \aref{small}). We would like to add that Claim 4.5 of \cite{MZajm10} and its proof given in Appendix D don't use the equation satisfied by $w$, so the proof is rigorously the same.\Box

\bigskip

In addition to the above lemmas, we have to add a new ingredient to the proof: the blow-up set of the multi-dimensional version $U(x,t)$ is always strictly under the forward light cone of vertex $(0,T(0))$.
More precisely, we make the following statement:
\begin{lem}\label{newlem}{\bf (The blow-up set is strictly under the forward light cone with vertex $(0,T(0))$)} For all $r>0$, we have $T(r)<T(0)+r$.
\end{lem}

\medskip

Let us first use the previous lemmas to derive Proposition \ref{pbut}, then we will prove Lemma \ref{newlem}.

\bigskip

{\it Proof of Proposition \ref{pbut} assuming that Lemma \ref{newlem} is true}:\\
(i) Arguing by contradiction, we assume that $\SS \cap \R^*_+$ contains an open non empty interval. Maximizing that interval, we can assume that a maximal interval $(a,b)$ is included in $\SS \cap \R^*_+$ with $0\le a< b \le +\infty$. From Lemma \ref{cns} and the fact that for all $r\ge 0$, $T(r)\ge 0$, we have two cases:

\medskip

\noindent{\bf Case 1}: $b=+\infty$ and for all $r\ge a$, $T(r) = T(a)+r-a$;\\
{\bf Case 2}: $0<b<+\infty$, $b\in\partial \SS$ and for all $r\in[a,b]$, $T(r) = T(c^*)-|r-c^*|$ for some $c^*\in[a,b]$.

\medskip

Let us find a contradiction in these two cases, assuming first that $a=0$ then $a>0$.

\medskip

- If $a=0$, then we see from Lemma \ref{newlem} that a contradiction follows in Case 1 or in Case 2 if $c^*>0$. Now, if Case 2 holds with $c^*=0$, then,
\begin{equation}\label{hassine}
\forall r\in[0,b],\;\;T(r) = T(0)-r\mbox{ and }b\in\R^*_+\cap\partial \SS.
\end{equation}
If $k(b)=0$, then a contradiction follows from (i) of Lemma \ref{lemfrank}.\\
If $k(b) =1$, then from the fact that $b\in\partial \SS$, there exists a sequence $r_n\in\RR$ converging to $b$. Since $(0,b)\in\SS$ by hypothesis, we have $r_n>b$ for $n$ large enough. Using (ii) of Lemma \ref{lemfrank}, we see that $b$ is right non-characteristic. Since $b$ is clearly left non-characteristic from \eqref{hassine}, this means that $b\in\RR\cap \R^*_+$. Since $\RR\cap \R^*_+$ is open from (i) of Theorem \ref{thbb}, this is in contradiction with the fact that $b\in\partial \SS$.\\
If $k(b) \ge 2$, then we know from Part 2 above that the upper bound in \eqref{corner} holds, which means that the blow-up set is corner shaped near $b$, and this is a contradiction by \eqref{hassine}.

\medskip

- If $a>0$, then, $a\in\partial\SS$ since the interval $(a,b)$ is maximal. If Case 2 holds with $c^*=a$, then the proof is exactly the same as in the case given above where $a=0$ and Case 4 holds with $c^*=0$. Now, if Case 1 holds or Case 2 with $c^*>a$ holds, 
then 
\[
\forall r\in[a,c^*],\;\;T(r) = T(a)+r-a\mbox{ and }a\in\R^*_+\cap\partial \SS,
\]
since the interval $(a,b)$ is maximal. The situation is symmetric with the situation where $a=0$ and Case 2 holds with $c^*=0$, and where we found a contradiction at the point $b$. Here, the contradiction follows in the same way, but at the point $a$. This concludes the proof of (i) in Proposition \ref{pbut}.

\bigskip

(ii) The proof is exactly the same as the proof of the analogous statement in \cite{MZajm10} (see The proof of Proposition 4.1 in Section 4.1 of that paper).\\
This concludes the proof of Proposition \ref{pbut} assuming that Lemma \ref{newlem} is true. It remains then to prove Lemma \ref{newlem}.\Box

\bigskip

{\it Proof of Lemma \ref{newlem}}: From invariance by time translation, we may assume that 
\[
T(0)=0.
\]
We proceed by contradiction and assume that for some $r_0>0$, we have $T(r_0)\ge r_0$. Recalling that $r\mapsto T(r)$ is $1$-Lipschitz, we see that $T(r_0)= r_0$. Using Lemma \ref{maxslope} (which is valid when $r_1=0$), we see that 
\begin{equation}\label{straight}
\forall r\in(0,r_0), T(r) = r\mbox{ and }(0,r_0)\subset \SS.
\end{equation}
In particular, by definition \eqref{nonchar} of a non-characteristic point,
\[
0\in\RR.
\]
Recalling that $p$ is subcritical or critical with respect to the existence of the conformal invariance (see \eqref{condp}), we have the following bound on the blow-up rate from Theorem 1 page 1149 in \cite{MZajm03} and Theorem 1 page 397 in \cite{MZma05}:
\begin{equation}\label{boundmulti}
\forall s\ge s_0,\;\;\frac 1K\le \|W_0(s)\|_{H^1(|Y|<1)}+\|\ps W_0(s)\|_{L^2(|Y|<1)}\le K,
\end{equation}
where $s_0\in\R$, $K>0$ and
\begin{equation}\label{wW}
\forall Y \in \R^N,\;\;W_0(Y,s) = w_0(|Y|,s).
\end{equation}
We claim that in order to conclude, it is enough to prove that 
\begin{equation}\label{hadaf}
\forall n\in\N,\;\;\int_{|Y|<1}|W_0(Y,s_n)|^{p+1}dY \ge \epsilon_0\mbox{ for some }s_n\to \infty
\end{equation}
and $\epsilon_0>0$.
Indeed, if \eqref{hadaf} holds, 
then we write from \eqref{boundmulti},
\begin{equation}\label{dalila}
\forall n\in\N,\;\;\int_\delta^{1-\delta}|w_0(y,s_n)|^{p+1}dy \ge \frac{\epsilon_0}2
\end{equation}
for some $\delta\in(0,\frac 12)$.\\
Following our argument for the proof of Claim 3.12 page 77 in \cite{MZcmp08}, we take $b\in(0,r_0)$. Recalling from \eqref{straight} that $T(b)=b$, we write by definition of similarity variables \eqref{defw} that
\begin{equation}\label{trans}
w_b(y',s') = (1-be^{s'})^{-\frac 2{p-1}}w_0(y,s)\mbox{ with }y=\frac{y'+be^{s'}}{1-be^{s'}}\mbox{ and }s=s'-\log(1-be^{s'}).
\end{equation}
Introducing 
\begin{equation}\label{defsn'}
s_n'=-\log(b+e^{-s_n}),\;\;y_1'(s')=\delta -b(1+\delta)e^{s'}\mbox{ and }y_2'(s')=1-\delta -b(2-\delta)e^{s'}
\end{equation}
we see that 
\begin{equation}\label{limit}
0<e^{s_n'}b<1\mbox{ and }e^{s_n'}b\to 1\mbox{ as }n\to\infty,
\end{equation} 
hence 
\[
y_i'(s_n')\in(-1,1),\;\;y_1(s) = \delta,\;\;y_2(s)=1-\delta,\;\;\forall y'\in(y_1'(s_n'), y_2(s_n')),\;\;\rho(y') \ge \frac{(1-be^{s'})^{\frac 2{p-1}}}{C(\delta)}.
\]
Therefore, we write from \eqref{trans} and \eqref{dalila}, 
\begin{eqnarray*}
&&\int_{-1}^1 |w_b(y',s_n')|^{p+1} \rho(y') dy' \ge \int_{y_1'(s_n')}^{y_2'(s_n')} |w_b(y',s_n')|^{p+1} \rho(y') dy'\\ 
&\ge& \frac{(1-be^{s_n'})^{-\frac {p+1}{p-1}}}{C(\delta)}\int_\delta^{1-\delta}|w_0(y,s_n)|^{p+1}dy\ge  \frac{\epsilon_0(1-be^{s_n'})^{-\frac {p+1}{p-1}}}{2C(\delta)}\to \infty,\mbox{ as }n\to \infty
\end{eqnarray*}
from \eqref{limit}.
This contradicts the bound
\[
\int_{-1}^1 |w_b(y',s_n')|^{p+1} \rho(y') dy' \le C_0(b)
\]
stated in Proposition \ref{boundedness}. Therefore, it is enough to prove \eqref{hadaf} in order to conclude the proof of Lemma \ref{newlem}.

\bigskip

Let us proceed by contradiction in order to prove \eqref{hadaf}, and assume that \begin{equation}\label{lp+1}
\int_{|Y|<1}|W_0(Y,s)|^{p+1}dY \to 0\mbox{ as }s\to \infty.
\end{equation}
Therefore,
 we see from \eqref{boundmulti} that for $s$ large enough, we have
\begin{equation}\label{updown}
\frac 1{2K}\le \|\nabla W_0(s)\|_{L^2(|Y|<1)}+\|\ps W_0(s)\|_{L^2(|Y|<1)}\le K.
\end{equation}
Introducing for all $n\in\N$,
\begin{equation}\label{defvn}
V_n(\xi,\tau) = (1-\tau)^{-\frac 2{p-1}}W_0(Y,s)\mbox{ with } Y = \frac \xi{1-\tau},\;\;s=n-\log(1-\tau),
\end{equation}
we see from the definitions \eqref{wW}, \eqref{defw} and \eqref{defu} of $W_0$, $w_0$ and $u$ that $V_n$ is radial in the sense that
\begin{equation}\label{vV}
V_n(\xi,\tau)=v_n(|\xi|,\tau)
\end{equation}
and that $V_n$ is a solution of the multi-dimensional equation \eqref{equrn} in the backward light cone of vertex $(0,1)$ above the section at time $\tau=0$, in the sense that
\[
\forall \tau \in[0,1),\;\;\forall |\xi|<1-\tau,\;\;\partial_\tau^2 V_n = \Delta V_n+|V_n|^{p-1}V_n.
\] 
Using \eqref{defvn}, \eqref{updown} and \eqref{lp+1}, we see that for $n$ large enough and for all $\tau \in[0,1)$,
\begin{equation}\label{l0}
\begin{array}{rcl}
\frac{(1-\tau)^{N-\frac{2(p+1)}{p-1}}}{128K^2}\le \|\partial_\tau V_n(\tau)\|_{L^2(|\xi|<1-\tau)}^2+\|\nabla V_n(\tau)\|_{L^2(|\xi|<1-\tau)}^2&\le& 16K^2 (1-\tau)^{N-\frac{2(p+1)}{p-1}},\\
\|V_n(\tau)\|_{L^{p+1}(|\xi|<1-\tau)}+\|V_n(\tau)\|_{L^2(|\xi|<1-\tau)}&\le & \epsilon_n(1-\tau)^{\frac N{p+1} - \frac 2{p-1}}
\end{array}
\end{equation}
where $\epsilon_n \to 0$ as $n\to \infty$.
Since $N-\frac{2(p+1)}{p-1}\le -2$ from the condition \eqref{condp} on $p$, there exists 
\begin{equation}\label{t0}
\tau_0(K)\in(0,1)
\end{equation}
 such that 
\begin{equation}\label{sommet}
\|\partial_\tau V_n(\tau_0)\|_{L^2(|\xi|<1-\tau_0)}^2+\|\nabla V_n(\tau_0)\|_{L^2(|\xi|<1-\tau_0)}^2\ge 20K^2.
\end{equation}
Applying the following one-dimensional Sobolev inequality 
\[
f(1)^2\le C \|f\|_{L^2(\frac 12,1)}\|f\|_{H^1(\frac 12,1)}
\]
to the radial version $v_n$ \eqref{vV}, we see from \eqref{t0} and \eqref{l0} that
\begin{equation}\label{lateral}
\sup_{(\xi,\tau)\in B_{\tau_0}}V_n(\xi,\tau)^2\le C(K)\sup_{0\le \tau \le \tau_0}\|V_n\|_{L^2(|\xi|<1-\tau)}\|V_n\|_{H^1(|\xi|<1-\tau)}\le C(K) \epsilon_n
\end{equation}
where 
\begin{equation}\label{defBt0}
B_{\tau_0}= \{(\xi,\tau)\;|\;0\le \tau \le \tau_0,\;|\xi|=1-\tau\}
\end{equation}
is the lateral boundary of the portion of the backward light cone of vertex $(0,1)$ located between the sections at $\tau=0$ and $\tau=\tau_0$.
Introducing the following local energy defined in the section of the backward light cone of vertex $(0,1)$
\[
\E(V_n(\tau)) =\int_{|\xi|<1-\tau}\left[(\partial_\tau V_n(\xi,\tau))^2+(\nabla V_n(\xi,\tau))^2-\frac{|V_n(\xi,\tau)|^{p+1}}{p+1}\right]d\xi,
\] 
we write from classical estimates, 
\[
\E(V_n(\tau_0)) -\E(V_n(0)) = -\frac 1{\sqrt 2}\int_{B_{\tau_0}}\left[\frac{|\nabla V_n - \frac \xi{|\xi|}\partial_\tau V_n|^2}2-\frac{|V_n|^{p+1}}{p+1}\right]d\sigma\le \int_{B_{\tau_0}}\frac{|V_n|^{p+1}}{\sqrt 2(p+1)}d\sigma
\]
where $B_{\tau_0}$ is defined in \eqref{defBt0}.
Using \eqref{l0}, \eqref{t0} and \eqref{lateral}, we see that for $n$ large enough, we have
\begin{eqnarray*}
&&\|\partial_\tau V_n(\tau_0)\|_{L^2(|\xi|<1-\tau_0)}^2+\|\nabla V_n(\tau_0)\|_{L^2(|\xi|<1-\tau_0)}^2\\
&\le& \|\partial_\tau V_n(0)\|_{L^2(|\xi|<1)}^2+\|\nabla V_n(0)\|_{L^2(|\xi|<1)}^2+K^2 \le 17K^2.
\end{eqnarray*}
This is a contradiction by \eqref{sommet}. Thus, \eqref{hadaf} holds and Lemma \ref{newlem} is proved. Since Proposition \ref{pbut} follows from Lemma \ref{newlem}, this concludes the proof of Lemma \ref{newlem} and Proposition \ref{pbut} too.\Box

\bigskip

{\bf Part 4: Proof of Theorem \ref{thgeo}}

The analogous statement in one space dimension is given in Theorems 1 and 2 
in \cite{MZisol10}. Thus, we need to say how to adapt the analysis of the paper \cite{MZisol10} to the radial case. Let us recall the strategy of the proof from Section 1.3 in that paper. Consider $u(r,t)$ a blow-up solution of equation \eqref{equ} and $r_0\in \SS\cap \R^*_+$. The decomposition of $w_{r_0}(y,s)$ is given in Theorem \ref{new} (up to replacing $u(r,t)$ by $-u(r,t)$, we may assume that $e^*_1=-1$).

\medskip
 
To prove that $r_0$ is an isolated characteristic point, the only tools we have are the energy criterion in (ii) of Corollary \ref{corcriterion} and the
trapping result stated in Theorem \ref{thtrap}. Note that due to the fact that equation \eqref{eqw} is no longer autonomous and depends on the considered blow-up point, the uniform version of the trapping we stated in this paper is strongly needed (we take $\rho_0=\frac{r_0}2$ in Theorem \ref{thtrap}).\\
In order to use these tools, we have to find the behavior of $w_r$ for $r$ near $r_0$. 
A simple idea is to start from the decomposition \eqref{cprofile00} for $w_{r_0}$ and the fact that the blow-up set is locally different from a straight line (it is in fact corner shaped by the upper bound in \eqref{corner} already proved in Parts 3 and 2 above), and use the transformation \eqref{defw} first to recover the behavior of $u(r,t)$, then the behavior of $w_r(y,s)$ for $r$ near $r_0$. Two problems arise in this simple idea:

\medskip

- we can't have information on $w_r(y,s)$ for all $y\in(-1,1)$, since information on the whole interval $(-1,1)$ would involve information on $w_{r_0}(y,s)$ for $|y|\ge 1$, and this is unavailable (at least at time $s$) because of the finite speed of propagation;

\medskip 

- the relation between $w_{r_0}$ and $w_r$ we get from \eqref{defw} depends explicitly on the value of $T(r)$ which is an unknown. The value of $T(r)$ specified by the upper bound in  \eqref{corner} changes the range of $s$ for which we have information.

\medskip

To overcome these problems, we proceed in 3 steps:

\medskip

{\bf Step 1: Initialization of the behavior of $w_r(s)$}

  Here, we use \eqref{cprofile00} and continuity arguments to show that for $r$ close enough to $r_0$ and $s=\lzero$ large enough, $w_r$ is close to a sum of $k$ solitons 
\[
\sum_{i=1}^k(-1)^i\kappa^*_1(\cd_i(\lzero),\tmu_i(\lzero))
\]
where 
\begin{equation}\label{defk*}
\kappa_1^*(d,\nu, y) =
\d\kappa_0\frac{(1-d^2)^{\frac 1{p-1}}}{(1+dy+\nu)^{\frac 2{p-1}}}.
\end{equation}
Note that $\kappa^*_1(d,\pm e^s)$ are heteroclinic orbits of the one-dimensional version of equation \eqref{eqw} connecting $\kappa(d)$ to $0$ or to $\infty$.

\medskip

{\bf Step 2: Propagation of the decomposition into solitons}

 Here, we are going to use essential facts of the theory of "solitons", namely
that under the flow of equation \eqref{eqw} and {\it uniformly} with respect to $r$ close to $r_0$, this decomposition is stable in time as time increases and that there are no collisions between solitons. More precisely, the idea is to use the equation \eqref{eqw} satisfied by $w_r$ to propagate this decomposition from $s=\lzero$ to $|\log|r-r_0||+L$ where $L$ is large, and prove (roughly speaking) the following (see Proposition 3.1 in \cite{MZisol10} for a precise statement): 
\begin{equation}\label{xdecomp}
\sup_{\lzero\le \ts \le |\log|r-r_0||+L}\left\|\vc{w_r(\ts)}{\ps w_r(s)}- \sum_{i=1}^{k}(-1)^{i}\kappa^*\left(\cd_i(s),\tmu_i(s)\right)\right\|_{\H}\to 0 
\end{equation}
as $L_{\k+1}\to \infty$, $L\to \infty$ and $r\to r_0$
for some parameters $(\cd_i(\ts),\tmu_i(\ts))$.\\
Let us remark that we can reduce to the case 
\[
r_0=T(r_0)=0,
\]
provided that we change equation \eqref{equ} by the following:
\[
\partial^2_{t} u =\partial^2_r u+\frac{(N-1)}{r_0+r} \pr u+|u|^{p-1}u\mbox{ for all }r>-r_0.
\]
This step involves two techniques:\\
- a modulation technique of the solution around the sum of solitons. The one-dimensional case is treated in Section 2 in \cite{MZisol10} which happens to be independent of the equation, hence it holds in the radial case with no modifications. More precisely, Section 2 depends only on the solitons \eqref{defk*} which are the same in one space dimension and in higher dimensions in the radial case.\\
- the proof of the stability of the decomposition into a decoupled sum of solitons, performed in Section 3 of \cite{MZisol10} in the one-dimensional case. If many estimates are independent from the equation, it is natural that some parts slightly change, because they use the equation in similarity variables. It is important to note that even though equation \eqref{eqw} depends on the considered point $r$ and on time $s$, the difference with the one-dimensional case comes only from one term whose effect is {\it uniformly} bounded thanks to \eqref{small}.
The estimate \eqref{small} makes it easy to control the effect of the additional term in the adaptation of Appendix C of \cite{MZisol10}, where we project the linearization of equation \eqref{eqw} around the sum of solitons. More precisely, as we did in Proposition \ref{lemproj}, we have to add the term $Ce^{-s}$ to the right-hand side of the four differential inequalities of Lemma C.2 in \cite{MZisol10}.
Accordingly, we have to mention that the statement of Claim 3.8 in \cite{MZisol10} slightly changes by adding the term $e^{-s}/\delta^*$ to the right-hand side of the differential inequalities involving $h_1$ and $h_2$, and also by mentioning that $\delta^*(L_m) \in (0,1)$ and not just that $\delta^*(L_m)>0$.

\medskip

{\bf Step 3: Trapping near one soliton and conclusion}

Following \aref{xdecomp}, it happens that at time $s=|\log|r-r_0||+L$, all the solitons $\kappa^*(\cd_i(s),\tmu_i(s))$ for $i=2,...,k$ become small (or vanish) for $L$ large and $|r-r_0|$ small (see Claim 3.4 in \cite{MZisol10} for a precise statement), 
so only the first soliton is left in \eqref{xdecomp}. 
%
Since 
\[
\forall r\in[\frac{r_0}2, \frac{3r_0}2]\mbox{ and }s\ge s_3(r_0),\;\;E(w_r(s))\ge E(\kappa_0) - C_3(r_0)e^{-s}
\]
from (i) of Corollary \ref{corcriterion},
we see that the first soliton has to be a pure soliton of the form $-\kappa(\cd^*_1,0)$ given in \eqref{defkd}, for some explicit $\cd^*_1=\cd_1^*(r)$, leading to the following estimate:
\[
\left\|\vc{w_r(s^*)}{\ps w_r(s^*)}+\vc{\kappa\left(\cd^*_1\right)}{0}\right\|_{\H}\le \epsilon_0\mbox{ where }s^*=|\log|r-r_0||+L
\]
and $\epsilon_0=\epsilon_0(\frac{r_0}2)>0$ is defined in the trapping result stated in (ii) of Theorem \ref{thtrap}. Here, there is a delicate point in the adaptation, since the statement of Theorem \ref{thtrap} is different from the one-dimensional case, in the sense that we need to apply it uniformly for $r\ge \rho_0\equiv \frac{r_0}2$, which is the case whenever $|r-r_0|$ is small. Applying that trapping result, we derive two facts:\\
- the point $r$ is non characteristic, hence $r_0$ is an isolated characteristic point (this is the conclusion of (i) in Theorem \ref{thgeo});\\
- the slope $T'(r)$ satisfies $|\argth T'(r) - \argth \cd_1^*(r)|\le C \epsilon_0$, which gives by integration the desired estimate in (ii) of Theorem \ref{thgeo}. 

\medskip

%
%
%
%

\def\cprime{$'$}

\noindent{\bf Address}:\\
Universit\'e de Cergy Pontoise, D\'epartement de math\'ematiques, 
2 avenue Adolphe Chauvin, BP 222, 95302 Cergy Pontoise cedex, France.\\
\vspace{-7mm}
\begin{verbatim}
e-mail: merle@math.u-cergy.fr
\end{verbatim}
Universit\'e Paris 13, Institut Galil\'ee, 
Laboratoire Analyse, G\'eom\'etrie et Applications, CNRS UMR 7539,
99 avenue J.B. Cl\'ement, 93430 Villetaneuse, France.\\
\vspace{-7mm}
\begin{verbatim}
e-mail: Hatem.Zaag@univ-paris13.fr
\end{verbatim}

\end{document}